\documentclass[12pt]{amsart}
\usepackage{graphicx}
\numberwithin{equation}{section} 

\newtheorem{thm}{Theorem}[section]
\newtheorem{cor}[thm]{Corollary}
\newtheorem{lem}[thm]{Lemma}

\theoremstyle{definition}

\theoremstyle{remark}

\newcommand{\bea}{\begin{eqnarray}}
\newcommand{\eea}{\end{eqnarray}}
\newcommand{\ba}{\begin{array}}
\newcommand{\ea}{\end{array}}
\newcommand{\bc}{\begin{center}}
\newcommand{\ec}{\end{center}}
\newcommand{\be}{\begin{equation}}
\newcommand{\ee}{\end{equation}}

\def\bn{{\mathbb N}}
\def\s{\sigma}

\def\l{\lambda}

\def\i{\varepsilon}
\def\t{\tau}
\def\xb{{\mathbf{x}}}
\def\cf{{\mathcal F}}
\def\cp{{\mathcal P}}

\def\br{\mathbb{R}}
\def\a{\alpha}

\def\m{\mu}

\begin{document}

\title[ON DOMINANT CONTRACTIONS]
{ON DOMINANT CONTRACTIONS AND A GENERALIZATION OF THE ZERO-TWO LAW}

\author{Farrukh Mukhamedov}
\address{Farrukh Mukhamedov\\
 Department of Computational \& Theoretical Sciences\\
Faculty of Science, International Islamic University Malaysia\\
P.O. Box, 141, 25710, Kuantan\\
Pahang, Malaysia} \email{{\tt far75m@yandex.ru} {\tt
farrukh\_m@iiu.edu.my}}

\begin{abstract}
Zaharopol proved the following result: let $T,S:L^1(X,{\cf},\m)\to
L^1(X,{\cf},\m)$ be two positive contractions such that $T\leq S$.
If $\|S-T\|<1$ then $\|S^n-T^n\|<1$ for all $n\in\bn$. In the
present paper we generalize this result to multi-parameter
contractions acting on $L^1$. As an application of that result we
prove a generalization of the "zero-two" law.

\vskip 0.3cm \noindent


\end{abstract}

\maketitle

\section{Introduction}

Let $(X,\cf,\m)$ be a measure space with a positive $\s$-additive
measure $\m$. In what follows for the sake of shortness by $L^1$ we
denote the usual  $L^1(X,\cf,\m)$ space associated with
$(X,\cf,\m)$. A linear operator $T:L^1\to L^1$ is called a {\it
positive contraction} if  $Tf\geq 0$ whenever $f\geq 0$ and
$\|T\|\leq 1$.

In \cite{OS} it was proved so called "zero-two" law for positive
contractions of $L^1$-spaces:

\begin{thm}\label{0-2} Let $T:L^1\to
L^1$ be a positive contraction.  If for some $m\in\bn\cup\{0\}$ one
has $\|T^{m+1}-T^m\|<2$, then
$$
\lim\limits_{n\to\infty}\|T^{n+1}-T^n\|=0.
$$
\end{thm}

In \cite{D} it was proved a "zero-two" law for Markov processes,
which allowed to study random walks on locally compact groups. Other
extensions and generalizations of the formulated law have been
investigated by many authors \cite{L1,F1,F2}.

Using certain properties of $L^1$-spaces Zaharopol \cite{Z} by means
of the following theorem reproved Theorem \ref{0-2}.

\begin{thm}\label{Z1} Let $T,S:L^1\to
L^1$ be two positive contractions such that $T\leq S$. If
$\|S-T\|<1$ then $\|S^n-T^n\|<1$ for all $n\in\bn$
\end{thm}

In the paper we provide an example (see Example 2) for which the
formulated theorem \ref{Z1} can not be applied. Therefore, we prove
a generalization of Theorem \ref{Z1} for multi-parameter
contractions acting on $L^1$. As a consequence of that result we
shall provide a generalization of the "zero-two" law. Similar
generalization has been considered in \cite{F2}.

\section{Dominant operators}

Let $T, S: L^{1}\rightarrow L^{1}$ be two positive contractions. We
write $T\leq S$ if $S-T$ is a positive operator. In this case we
have
\begin{equation}\label{3.2}
\|Sx-Tx\|=\|Sx\|-\|Tx\|,
\end{equation}
for every $x\geq 0$.  Moreover, for positive operator $T:
L^{1}\rightarrow L^{1}$ one can prove the following equality
\begin{equation}\label{3.1}\|T\|=\sup_{\|x\|=1}\|Tx\|=\sup_{\|x\|=1,x\geq
0}\|Tx\|. \end{equation}

The main result of this section is the following

\begin{thm}\label{ZN0} Let
$T_1,T_2,S_1,S_2:L^{1}\to L^{1}$ be positive contractions such that
$T_i\leq S_i$, $i=1,2$ and $S_1S_2=S_2S_1$. If there is an
$n_{0}\in\mathbb{N}$ such that
$\|S_1S_2^{n_{0}}-T_1T_2^{n_{0}}\|<1$. Then
$\|S_1S_2^{n}-T_1T_2^{n}\|<1$ for every $n\geq n_{0}.$
\end{thm}

\begin{proof} Let us assume that $\|S_1S_2^{n}-T_1T_2^{n}\|=1$ for some $n>n_{0}.$
Therefore, denote
\begin{equation*}
m=\min\{n\in\mathbb{N}:\|S_1S_2^{n_{0}+n}-T_1T_2^{n_{0}+n}\|=1\}.
\end{equation*}
It is clear that $m\geq 1$. The inequalities $T_1\leq S_1$, $T_2\leq
S_2$ imply that $S_1S_2^{n_{0}+n}-T_1T_2^{n_{0}+n}$ is a positive
operator. Then according to \eqref{3.1} there exists a sequence
$\{x_{n}\}\in L^{1}$ such that $x_{n}\geq 0$, $\|x_{n}\|=1, \forall
n\in\mathbb{N}$ and
\begin{eqnarray}\label{eq1}
\lim\limits_{n\to\infty}\|(S_1S_2^{n_{0}+n}-T_1T_2^{n_{0}+n})x_{n}\|=1.
\end{eqnarray}

Positivity of $S_1S_2^{n_{0}+n}-T_1T_2^{n_{0}+n}$ and $x_{n}\geq 0$
together with  \eqref{3.2} imply that
\begin{eqnarray}\label{eq2}
\|(S_1S_2^{n_{0}+n}-T_1T_2^{n_{0}+n})x_{n}\|=\|S_1S_2^{n_{0}+m}x_{n}\|-\|T_1T_2^{n_{0}+m}x_{n}\|
\end{eqnarray}
for every $n\in\bn$. It then follows from \eqref{eq1},\eqref{eq2}
that

\begin{eqnarray}\label{eq3}
&&\lim\limits_{n\to\infty}\|S_1S_2^{n_{0}+m}x_{n}\|=1,\\
\label{eq4} &&\lim\limits_{n\to\infty}\|T_1T_2^{n_{0}+m}x_{n}\|=0.
\end{eqnarray}

Thanks to the contractivity of $S$, $Z$ and $S_1S_2=S_2S_1$ one gets
\begin{equation*}
\|S_1S_2^{n_0+m}x_n\|=\|S_2(S_1S_2^{n_0+m-1}x_n)\|\leq\|S_1S_2^{n_0+m-1}x_n\|\leq
\|S_2^{m}x_n\|
\end{equation*}
which with \eqref{eq3} yields
\begin{equation}\label{ZS1}
\lim\limits_{n\to\infty}\|S_1S_2^{n_0+m-1}x_n\|=1, \ \ \
\lim\limits_{n\to\infty}\|S_2^{m}x_n\|=1.
\end{equation}

Moreover, the contractivity of $S_i$, $T_i$ ($i=1,2$) implies that
$\|T_1T_2^{n_{0}+m-1}x_{n}\|\leq 1$, $\|T_2^{m}x_{n}\|\leq 1$ and
$\|S_1S_2^{n_{0}}T^{m}x_{n}\|\leq 1$ for every $n\in\mathbb{N}$.
Therefore, we may choose a subsequence $\{y_{k}\}$ of $\{x_{n}\}$
such that the sequences $\{\|T_1T_2^{n_{0}+m-1}y_{k}\|\}$,
$\{\|T_2^{m}y_k\|\}$, $\{\|S_1S_2^{n_{0}}T^{m}x_k\|\}$ converge. Put
\begin{eqnarray}\label{eq5}
&&\alpha=\lim\limits_{k\to\infty}\|T_1T_2^{n_{0}+m-1}y_{k}\|,\\
\label{eq6}
&&\beta=\lim\limits_{k\to\infty}\|S_1S_2^{n_{0}}T^{m}y_{k}\|,\\
\label{eq7} &&\gamma=\lim\limits_{k\to\infty}\|T_2^{m}y_{k}\|.
\end{eqnarray}

The inequality $\|S_1S_2^{n_{0}+m-1}-T_1T_2^{n_{0}+m-1}\|<1$ with
\eqref{ZS1} implies that $\alpha>0$. Hence  we may choose a
subsequence $\{z_{k}\}$ of $\{y_{k}\}$ such that
$\|T_1T_2^{n_{0}+m-1}z_{k}\|\neq 0$ for all $k\in\mathbb{N}$.

From $\|T_1T_2^{n_{0}+m-1}z_{k}\|\leq \|T_2^{m}z_{k}\|$ together
with \eqref{eq5}, \eqref{eq7} we find $\alpha\leq\gamma$, and hence
$\gamma>0.$

Using \eqref{3.2} one gets
\begin{eqnarray}\label{eq10}
\|S_1S_2^{n_{0}}T_2^{m}z_{k}\|&=&\|S_1S_2^{n_{0}+m}z_{k}-(S_1S_2^{n_{0}+m}z_{k}-S_1S_2^{n_{0}}T_2^{m}z_{k})\|\nonumber\\
&=&\|S_1S_2^{n_{0}+m}z_{k}\|-\|S_1S_2^{n_{0}+m}z_{k}-S_1S_2^{n_{0}}T_2^{m}z_{k}\|\nonumber\\
&\geq&
\|S_1S_2^{n_{0}+m}z_{k}\|-\|S_2^{m}z_{k}-T_2^{m}z_{k}\|\nonumber\\
&=&\|S_1S_2^{n_{0}+m}z_{k}\|-\|S_2^{m}z_{k}\|+\|T_2^{m}z_{k}\|
\end{eqnarray}
Due to \eqref{eq3},\eqref{ZS1} we have
$$\lim\limits_{k\to\infty}\|S_1S_2^{n_{0}+m}z_{k}\|-\|S_2^{m}z_{k}\|=0;$$
which with \eqref{eq10} implies that
$$\lim\limits_{k\to\infty}\|S_1ZS_2^{n_{0}}T_2^{m}z_{k}\|\geq\lim\limits_{k\to\infty}\|T_2^{m}z_{k}\|,$$
therefore, $\beta\geq\gamma.$

On the other hand, by
$\|S_1S_2^{n_{0}}T_2^{m}z_{k}\|\leq\|T_2^{m}z_{k}\|$ one gets
$\gamma\geq\beta$, hence $\gamma=\beta$.

Now  set
$$u_{k}=\frac{T_2^{m}z_{k}}{\|T_2^{m}z_{k}\|}, \ \ k\in\bn.$$
Then using the equality $\gamma=\beta$ and \eqref{eq4}  one has
\begin{eqnarray*}
&&\lim\limits_{k\to\infty}\|S_1S_2^{n_{0}}u_{k}\|=
\lim\limits_{k\to\infty}\frac{\|S_1S_2^{n_{0}}T_2^{m}z_{k}\|}{\|T_2^{m}z_{k}\|}=1,\\
&&\lim\limits_{k\to\infty}\|T_1T_2^{n_{0}}u_{k}\|=
\lim\limits_{k\to\infty}\frac{\|T_1T^{n_{0}+m}z_{k}\|}{\|T_2^{m}z_{k}\|}=0.
\end{eqnarray*}

So, owing to \eqref{3.2} and positivity of
$S_1S_2^{n_{0}}-T_1T_2^{n_{0}}$, we get

$$\lim\limits_{k\to\infty}\|(S_1S_2^{n_{0}}-T_1T_2^{n_{0}})z_{k}\|=1.$$
Since $\|u_{k}\|=1, u_{k}\geq 0, \forall k\in\mathbb{N}$ from
\eqref{3.1} one finds $\|S_1S_2^{n_{0}}-T_1T_2^{n_{0}}\|=1,$ which
is a contradiction. This completes the proof.
\end{proof}

\begin{cor}\label{ZN} Let
$Z,T,S:L^{1}\to L^{1}$ be positive contractions such that $T\leq S$
and $ZS=SZ$. If there is an $n_{0}\in\mathbb{N}$ such that
$\|Z(S^{n_{0}}-T^{n_{0}})\|<1$. Then $\|Z(S^{n}-T^{n})\|<1$ for
every $n\geq n_{0}.$
\end{cor}

Assume that $Z=Id$. If $n_0=1$, then from Corollary \ref{ZN} we
immediately get the Zaharopol's result (see Theorem \ref{Z1}).  If
$n_0>1$ then we obtain a main result of \cite{MHT}.

Let us provide an example of $Z,S,T$ positive contractions for which
statement of Corollary \ref{ZN} is satisfied.

{\bf Example 1.} Consider $\br^2$ with a norm $\|\xb\|=|x_1|+|x_2|$,
where $\xb=(x_1,x_2)$. An order in $\br^2$ is defined as usual,
namely $\xb\geq 0$ if and only if $x_1\geq 0$, $x_2\geq 0$. Now
define mappings $Z:\br^2\to\br^2$,$T:\br^2\to\br^2$ and
$S:\br^2\to\br^2$, respectively, by
\begin{eqnarray}\label{Z}
&& Z(x_1,x_2)=(ux_1+vx_2,ux_2),\\
\label{S}
&& S(x_1,x_2)=\bigg(\frac{x_1+x_2}{2},\frac{x_2}{2}\bigg),\\
\label{T} && T(x_1,x_2)=(\l x_2,0).
\end{eqnarray}
 The positivity of $Z$,$S$ and $T$ implies that $u,v,\l\geq 0$. It
is easy to check that $T\leq S$ holds if and only if $2\l\leq 1$.

One can see that
\begin{eqnarray*}
\|Z\|=\sup_{\|\xb\|=1\atop \xb\geq 0}\|Z\xb\|&=&\max_{x_1+x_2=1\atop
x_1,x_2\geq 0}
\{ux_1+(u+v)x_2\}\nonumber\\
&=&\max_{0\leq x_2\leq 1}\{u+vx_2\}\nonumber\\
&=&u+v
\end{eqnarray*}
Hence, contractivity of $Z$ implies that $u+v=1$. Similarly, we find
that $\|S\|=1$ and $\|T\|=\l$. From \eqref{Z} and \eqref{S} one gets
that $ZS=SZ$.

By means of \eqref{Z},\eqref{S},\eqref{T} one finds Similarly, one
gets
\begin{eqnarray}\label{S-T}
\|Z(S-T)\|=\sup_{\|\xb\|=1\atop \xb\geq 0}\|Z(S-T)\xb\|&=&
\max_{x_1+x_2=1\atop x_1,x_2\geq 0}\bigg\{\frac{1}{2}\big(ux_1+x_2+ux_2-2\l ux_2\big)\bigg\}\nonumber\\
&=&\frac{1+u(1-2\l)}{2}.
\end{eqnarray}
The condition $2\l\leq 1$ yields that $\|Z(S-T)\|<1$. Consequently,
Corollary \ref{ZN} implies $\|Z(S^n-T^n)\|<1$ for all $n\in\bn$.

Now let us formulate a multi-parametric version of Theorem
\ref{0-2}.

\begin{thm}\label{ZN-m} Let
$T_i,S_i:L^{1}\to L^{1}$, $i=1,\dots,N$ be positive contractions
such that $T_i\leq S_i$ with
\begin{equation}\label{zn-c}T_iT_j=T_jT_i, \quad
S_iS_j=S_jS_i \ \ \textrm{for every} \ i,j=1,\dots,N.
\end{equation}
If there are $n_{i,0}\in\mathbb{N}$, $i=1,\dots,N$ such that
\begin{equation}\label{zn-1}
\|S_1^{n_{1,0}}\cdots S_N^{n_{N,0}}-T_1^{n_{1,0}}\cdots
T_N^{n_{N,0}}\|<1.
\end{equation}
 Then
\begin{equation}\label{zn-2}
 \|S_1^{m_1}\cdots
S_N^{m_{N}}-T_1^{m_1}\cdots T_N^{m_{N}}\|<1
\end{equation}
for all
$m_i\geq n_{i,0}$, $i=1,\dots,N$.
\end{thm}

\begin{proof} Let us fix the first $N-1$ operators in \eqref{zn-1},
i.e. for a moment we denote
\begin{equation}\label{zn-3}
\mathbf{S}_{N-1}=S_1^{n_{1,0}}\cdots S_{N-1}^{n_{N-1,0}} \quad
\mathbf{T}_{N-1}=T_1^{n_{1,0}}\cdots T_{N-1}^{n_{N-1,0}},
\end{equation}
then \eqref{zn-1} can be written as follows
\begin{equation*}\label{zn-4}
\|\mathbf{S}_{N-1}S_N^{n_{N,0}}-\mathbf{T}_{N-1}T_N^{n_{N,0}}\|<1.
\end{equation*}
After applying Theorem \ref{ZN0} to the last inequality we find
\begin{equation}\label{zn-4}
\|\mathbf{S}_{N-1}S_N^{m_N}-\mathbf{T}_{N-1}T_N^{m_{N}}\|<1
\end{equation}
for all $m_N\geq n_{N,0}$. Now taking into account \eqref{zn-3} and
\eqref{zn-c} we rewrite \eqref{zn-4} as follows
\begin{equation}\label{zn-5}
\|S_N^{m_{N}}S_1^{n_{1,0}}\cdots
S_{N-1}^{n_{N-1,0}}-T_N^{m_{N}}T_1^{n_{1,0}}\cdots
T_{N-1}^{n_{N-1,0}}\|<1.
\end{equation}
Now again applying the same idea as above to \eqref{zn-5} we get
\begin{equation*}
\|S_N^{m_{N}}S_1^{n_{1,0}}\cdots
S_{N-1}^{m_{N-1}}-T_N^{m_{N}}T_1^{n_{1,0}}\cdots
T_{N-1}^{m_{N-1}}\|<1,
\end{equation*}
for all $m_{N-1}\geq n_{N-1,0},m_N\geq n_{N,0}$. Hence, continuing
this procedure $N-2$ times we obtain the desired inequality.
\end{proof}

{\bf Remark 3.1.}  It should be noted the following:
\begin{enumerate}
\item[(i)] Since the dual of $L^1$ is $L^\infty$ then due to the
duality theory the proved Theorems \ref{ZN0} and \ref{ZN-m} holds
true if we replace $L^1$-space with $L^\infty$.

\item[(ii)] Unfortunately, that the proved theorems and its corollaries are not
longer true if one replaces $L^1$-space by an $L^p$-space,
$1<p<\infty$. Indeed, consider $X=\{1,2\}$, $\cf=\cp(\{1,2\})$ and
the measure $\m$ is given by $\m(\{1\})=\m(\{2\})=1/2$. In this
case, $L^p$ is isomorphic to the Banach lattice $\br^2$ (here an
order is defined as usual, namely $\xb\geq 0$ if and only if
$x_1\geq 0$, $x_2\geq 0$) with the norm
$\|\xb\|_p=\big(|x_1|^p+|x_2|^p\big)^{1/p}/2$, where
$\xb=(x_1,x_2)$. Define two operators by
\begin{eqnarray*}
&& S(x_1,x_2)=\bigg(\frac{x_1+x_2}{2},\frac{x_1+x_2}{2}\bigg), \quad
T(x_1,x_2)=\bigg(0,\frac{x_1}{2}\bigg)
\end{eqnarray*}
Then it is shown (see \cite{Z}) that $\|S-T\|<1$, but
$\|S^2-T^2\|=1$.

\item[(iii)] It would be better to note that certain ergodic properties
of dominant positive operators has been studied in \cite{E}. In
general, a monograph \cite{K} is devoted to dominant operators.

\end{enumerate}

Let us give another example, for which conditions of Theorem
\ref{Z1} does not hold, but Theorem \ref{ZN0} can be applied.

{\bf Example 2.} Let us consider $\br^2$ as in Example 1.  Now
define mappings $T:\br^2\to\br^2$ and $S:\br^2\to\br^2$ as follows
\begin{eqnarray}\label{S1}
&& S(x_1,x_2)=\bigg(\frac{1}{2}x_1+\frac{1}{3}x_2,\frac{1}{2}x_1+\frac{1}{3}x_2\bigg),\\
\label{T1} && T(x_1,x_2)=\bigg(\frac{1}{4}x_2,0\bigg).
\end{eqnarray}
It is clear that $S$ and $T$ are positive and $T\leq S$.

One can see that $\|S\|=1$, $\|T\|=1/4$. From \eqref{S1},\eqref{T1}
one gets
\begin{eqnarray}\label{S-T}
&&\|S-T\|=\sup_{\|\xb\|=1\atop \xb\geq 0}\|(S-T)\xb\|= \max_{0\leq
x_1 \leq 1}\bigg\{\frac{7x_1+5}{12}\bigg\}=1 \\[2mm]\label{S-T2}
&&\|S^2-T^2\|=\sup_{\|\xb\|=1\atop \xb\geq 0}\|(S^2-T^2)\xb\|=
\max_{0\leq x_1\leq 1}\bigg\{\frac{5x_1+10}{18}\bigg\}=\frac{15}{18}
\end{eqnarray}

Consequently, we have positive contractions $T$ and $S$ with $S\geq
T$ such that $\|S-T\|=1$,$\|S^2-T^2\|<1$. This shows that the
condition of Theorem \ref{Z1} is not satisfied, but due to Corollary
\ref{ZN} with $Z=id$ we have $\|S^n-T^n\|<1$ for all $n\geq 2$.
Therefore the
proved Theorem \ref{ZN} is an extension of the Zaharopol's result.\\

\section{A generalization of the zero-two law}

In this section we are going to prove a generalization of the
zero-two law for positive contractions on $L^1$. Before formulate
the main result we prove some auxiliary facts.

First note that for any $x,y\in L^1$ one defines
\begin{equation}\label{wedge}
x\wedge y=\frac{1}{2}(x+y-|x-y|).
\end{equation}
It is well known (see \cite{AB}) that for any mapping $S$ of $L^1$
one can define its modulus by
\begin{equation}\label{modul}
|S|x=\sup\{ Sy : \ |y|\leq x\}, \ \ x\in L^1, x\geq 0.
\end{equation}
Hence, similarly to \eqref{wedge} for given two mappings $S,T$ of
$L^1$ we define
\begin{equation}\label{wedge2}
(S\wedge T)x=\frac{1}{2}(Sx+Tx-|S-T|x), \ \ x\in L^1.
\end{equation}

A linear operator $Z:L^1\to L^1$ is called {\it a lattice
homomorphism} whenever
\begin{equation}\label{V}
Z(x\vee y)=Zx\vee Zy
\end{equation}
holds for all $x,y\in L^1$. One can see that such an operator is
positive. Note that such homomorphisms were studied in \cite{AB}.

Recall that a net $\{x_\a\}$ in $L^1$ is {\it order convergent} to
$x$, denoted $x_\a\to ^o x$ whenever there exists another net
$\{y_\a\}$ with the same index set satisfying $|x_\a - x|\leq
y_\a\downarrow 0$. An operator $T:L^1\to L^1$ is said to be {\it
order continuous}, if $x_\a\to ^o 0$ implies $Tx_\a\to^o 0$.

\begin{lem}\label{1modul}
Let $S$, $T$ be positive contractions of $L^1$, and $Z$ be an order
continuous lattice homomorphism of $L^1$. Then one has
\begin{eqnarray}\label{modul2}
Z|S-T|=|Z(S-T)|.
\end{eqnarray}
Moreover, we have
\begin{eqnarray}\label{wedge3}
Z(S\wedge T)=ZS\wedge ZT.
\end{eqnarray}
\end{lem}

\begin{proof} From \eqref{modul} we find that
\begin{eqnarray}\label{modul2}
Z|S-T|x&=&Z(\sup\{ (S-T)y : \ |y|\leq x\})\nonumber\\
&=&\sup\{ Z(S-T)y : \ |y|\leq x\})\nonumber\\
&=&|Z(S-T)|x,
\end{eqnarray}
for every $ x\in L^1, x\geq 0$.

The equality \eqref{wedge2} yields that
\begin{eqnarray}
&&Z|S-T|=ZS+ZT-2Z(S\wedge T),\nonumber\\
\label{modul4} &&|Z(S-T)|=ZS+ZT-2(ZS\wedge ZT),
\end{eqnarray}
which with \eqref{modul2} imply that
\begin{eqnarray*}
Z(S\wedge T)=ZS\wedge ZT.
\end{eqnarray*}

\end{proof}

In what follows, an order continuous lattice homomorphism $Z:L^1\to
L^1$ with $\|Z\|\leq 1$, is called {\it a lattice contraction}.

 Now we have the following

\begin{lem}\label{1-2} Let $Z$ be a lattice contraction and $T$ be a positive
contraction of $L^1$ such that $ZT=TZ$. If for some
$m\in\bn\cup\{0\}$, $k\in\bn$ one has $\|Z(T^{m+k}-T^m)\|<2$, then
$\|Z(T^{m+k}-T^{m+k}\wedge T^m)\|<1$.
\end{lem}

\begin{proof} According to the assumption there is $\delta>0$ such that
$\|Z(T^{m+k}-T^m)\|=2(1-\delta)$. Let us suppose that
$\|Z(T^{m+k}-T^{m+k}\wedge T^m)\|=1$. Then thanks to \eqref{3.1}
there exists $x\in L^1$ with $x\geq 0$, $\|x\|=1$ such that
$$
\|Z(T^{m+k}-T^{m+k}\wedge T^m)x\|>1-\frac{\delta}{4},
$$
which with \eqref{3.2} implies that $\|ZT^{m+k}x\|>1-\delta/4$ and
$\|Z(T^{m+k}\wedge T^m)x\|<\delta/4$.  The commutativity $T$ and $Z$
yields that $\|ZT^{m}x\|>1-\delta/4$.

Now using \eqref{modul4} and \eqref{wedge3}  one finds
\begin{eqnarray*}
\big\||Z(T^{m+k}-T^m)|x\big\|&=&\|ZT^{m+k}x\|+\|ZT^mx\|-2\|Z(T^{m+k}\wedge T^m)x\|\\
&>&1-\frac{\delta}{4}+1-\frac{\delta}{4}-2\cdot\frac{\delta}{4}\\
&=&2\bigg(1-\frac{\delta}{2}\bigg).
\end{eqnarray*}
This with the equality
$$
\big\||Z(T^{m+k}-T^m)|\big\|=\|Z(T^{m+k}-T^m)\|,
$$
contradicts to $\|Z(T^{m+k}-T^m)\|=2(1-\delta/2)$.
\end{proof}

\begin{lem}\label{1-3} Let $Z$ be a lattice contraction and $T$ be a positive contraction
of $L^1$ such that $ZT=TZ$. If for some $m\in\bn\cup\{0\}$,
$k\in\bn$ one has $\|Z(T^{m+k}-T^{m+k}\wedge T^m)\|<1$, then for any
$\i>0$ there are $d,n_0\in\bn$ such that
$$
\|Z^d(T^{n+k}-T^n)\|<\i \ \ \ \ \textrm{for all} \ n\geq n_0
$$
\end{lem}

\begin{proof} It is known that (see \cite{Z1}, p. 310)  for any contraction $T$ on
$L^1$ there is $\gamma>0$ such that
\begin{equation}\label{G-l}
\bigg\|\bigg(\frac{I+T}{2}\bigg)^\ell-T\bigg(\frac{I+T}{2}\bigg)^\ell\bigg\|
\leq\frac{\gamma}{\sqrt{\ell}}.
\end{equation}

Then for given $k\in \bn$, using \eqref{G-l} one easily finds that
\begin{equation}\label{G-l2}
\bigg\|\bigg(\frac{I+T}{2}\bigg)^\ell-T^k\bigg(\frac{I+T}{2}\bigg)^\ell\bigg\|
\leq\frac{k\gamma}{\sqrt{\ell}}.
\end{equation}

Let $\i>0$ and fix $\ell\in\bn$ such that
$k\gamma/\sqrt{\ell}<\i/4$.

Then according to Corollary \ref{ZN} from the assumption of the
lemma we have
\begin{equation}\label{Zl}
\big\|Z(T^{\ell(m+k)}-(T^{m+k}\wedge T^m)^\ell)\big\|<1.
\end{equation}

Hence,
\begin{eqnarray}\label{Zl2}
\bigg\|Z\bigg(T^{\ell(m+k)}&-&\bigg(\frac{I+T}{2}\bigg)^\ell(T^{m+k}\wedge
T^m)^\ell\bigg)\bigg\|=\nonumber\\&=&
\bigg\|Z\bigg(T^{\ell(m+k)}-\frac{1}{2^\ell}\sum_{i=0}^{\ell}C^i_\ell
T^i(T^{m+k}\wedge T^m)^\ell\bigg)\bigg\| \nonumber\\
&\leq&\sum_{i=0}^{\ell}\frac{C^i_\ell}{2^\ell}\big\|Z(T^{\ell(m+k)}-T^i(T^{m+k}\wedge T^m)^\ell)\big\| \nonumber\\
&\leq& \frac{1}{2^\ell}\big\|Z(T^{\ell(m+k)}-(T^{m+k}\wedge
T^m)^\ell)\big\|+\sum_{i=0}^{\ell}\frac{C^i_\ell}{2^\ell}\nonumber
\\
&<& \frac{1}{2^\ell}+\sum_{i=1}^{\ell}\frac{C^i_\ell}{2^\ell}=1.
\end{eqnarray}

Define
$$
Q_\ell:=T^{\ell(m+k)}-\bigg(\frac{I+T}{2}\bigg)^\ell(T^{m+k}\wedge
T^m)^\ell
$$
and put $V_\ell^{(1)}=(T^{m+k}\wedge T^m)^\ell$. Then one can see
that
$$
T^{\ell(m+k)}=\bigg(\frac{I+T}{2}\bigg)^\ell V_\ell^{(1)}+Q_\ell.
$$

Now for every $d\in\bn$, define
$$
V_\ell^{(d+1)}=T^{\ell(m+k)}V_\ell^{(d)}+V_\ell^{(1)}Q^d_\ell.
$$

Then by induction one can establish \cite{Z1} that
\begin{equation}\label{Td}
T^{d\ell(m+k)}=\bigg(\frac{I+T}{2}\bigg)^\ell V_\ell^{(d)}+Q^d_\ell
\end{equation}
for every $d\in\bn$.

Due to Proposition 2.1 \cite{Z} one has
\begin{equation}\label{V}
\|V^{(d)}_\ell\|\leq 2
\end{equation}
for all $d\in\bn$.

Now from \eqref{Zl2} we find $\|ZQ_\ell\|<1$, therefore there exists
$d\in \bn$ such that $\|(ZQ_\ell)^d\|<\i/4$. So, commutativity $Z$
and $T$ implies that $ZQ_\ell=Q_\ell Z$, which yields that $\|Z^d
Q_\ell^d\|<\i/4$.

Put $n_0=d\ell(m+k)$, then from \eqref{Td} with
\eqref{G-l2},\eqref{V} we get
\begin{eqnarray*}
\|Z^d(T^{n_0+k}-T^{n_0})\|&=&\bigg\|Z^d\bigg(T^k\bigg(\frac{I+T}{2}\bigg)^\ell-
\bigg(\frac{I+T}{2}\bigg)^\ell\bigg)V^{(d)}_\ell\\
&&+Z^d(T^kQ_\ell^d-Q^d_\ell)\bigg\|\\
&\leq&\bigg\|\bigg(T^k\bigg(\frac{I+T}{2}\bigg)^\ell-
\bigg(\frac{I+T}{2}\bigg)^\ell\bigg)V^{(d)}_\ell\bigg\|\\&&+\|Z^d
Q_\ell^d(T-1)\|\\
&\leq& 2\cdot\frac{k\gamma}{\sqrt{\ell}}+2\cdot\frac{\i}{4}<\i.
\end{eqnarray*}

Take any $n\geq n_0$, then from the last inequality one finds
\begin{equation*}
\|Z^d(T^{n+k}-T^n)\|=\|T^{n-n_0}Z^d(T^{n_0+k}-T^{n_0})\|\leq
\|Z^d(T^{n_0+k}-T^{n_0})\|<\i
\end{equation*}
which completes the proof.
\end{proof}

Now we are ready to formulate the main result of this section.

\begin{thm}\label{0-2-1} Let $Z$, $T$ be two positive contractions
of $L^1$ such that $TZ=ZT$. If for some $m\in\bn\cup\{0\}$,
$k\in\bn$ one has $\|Z(T^{m+k}-T^m)\|<2$, then for any $\i>0$ there
are $d,n_0\in\bn$ such that
$$
\|Z^d(T^{n+k}-T^n)\|<\i \ \ \ \ \textrm{for all} \ n\geq n_0
$$
\end{thm}

The proof of this theorem immediately follows from Lemmas \ref{1-2}
and \ref{1-3}.

 {\bf Remark.} Note that if we take as $Z=I$,$k=1$
then we obtain Theorem \ref{0-2} as a corollary of Theorem
\ref{0-2-1}.


\end{document}